\input amssym.def
\input amssym
\magnification=1200
\parindent0pt
\hsize=16 true cm
\baselineskip=13  pt plus .2pt
$ $

\def\G{(\Gamma,{\cal G})}
\def\Z{{\Bbb Z}}
\def\H{{\Bbb H}}
\def \S {{\cal S}}

\centerline {\bf On finite groups of isometries of handlebodies in arbitrary dimensions}

\centerline {\bf and finite extensions of Schottky groups}

\medskip

\centerline {by}

\medskip

\centerline {{\bf Mattia Mecchia}  and  {\bf Bruno P. Zimmermann}}

\bigskip \bigskip

{\bf Abstract.}  {\sl  It is known that  the order of a finite group of
diffeomorphisms of a 3-dimensional handlebody of genus $g>1$ is bounded by the
linear polynomial $12(g-1)$, and that the order of a finite group of
diffeomorphisms of a 4-dimensional handlebody (or equivalently, of its boundary
3-manifold), faithful on the fundamental group, is bounded by a quadratic
polynomial in $g$ (but not by a linear one).  In the present paper we prove a
generalization for handlebodies of arbitrary dimension $d$, uniformizing
handlebodies by Schottky groups and considering finite groups of isometries of
such handlebodies. We prove that the order of a finite group of isometries of a
handlebody of dimension $d$ acting faithfully on the fundamental group is
bounded by a polynomial of degree $d/2$ in $g$ if $d$ is even, and of degree
$(d+1)/2$ if $d$ is odd, and that the degree $d/2$ for even $d$ is best
possible.  This implies then analogous polynomial Jordan-type bounds for
arbitrary finite groups of isometries of handlebodies (since a handlebody of
dimension $d>3$ admits $S^1$-actions, there does not exist an upper bound for
the order of the group itself).}

\bigskip

{\sl 2010 Mathematics Subject Classification:}  57S17,  57S25, 57N16

{\sl Key words and phrases:}  handlebody, finite group action, Schottky group, Jordan-type
bound.

\bigskip

{\bf 1. Introduction.}  All finite group actions in the present paper will be faithful,
smooth and orientation-preserving, and all manifolds  will be orientable. We study finite
group actions of large order on handlebodies of dimension $d \ge 3$ and genus
$g > 1$.

\medskip

An orientable handlebody $V^d_g$ of dimension $d$ and genus $g$ can be defined
as a regular neighbourhood of a finite graph with free fundamental group of
rank $g$ embedded in the sphere $S^d$; alternatively, it is obtained from the
ball $B^d$ by attaching along its boundary  $g$ copies of a handle $B^{d-1}
\times [0,1]$ in an orientable way, or as the boundary-connected sum of $g$
copies of $B^{d-1} \times S^1$. The boundary of $V^d_g$ is a closed manifold
$H^{d-1}_g$ which is the connected sum of $g$ copies of $S^{d-2} \times S^1$.

\medskip

By [Z1] the order of a finite group of diffeomorphisms of a 3-dimensional
handlebody $V^3_g$ of genus $g>1$ is bounded by the linear polynomial $12(g-1)$
(see also [MMZ, Theorem 7.2], [MZ]); also, a finite group  $G$ acting
faithfully on $V^3_g$ acts faithfully also on the fundamental group.  On the
other hand, since the closed 3-manifold $H^3_g$ admits $S^1$-actions, it has
finite cyclic group actions of arbitrarily large order acting trivially on the
fundamental group, and the same is true also for handlebodies $V^d_g$ of
dimensions $d > 3$. However it is shown in [Z4] that if a finite group of
diffeomorphisms of $H^3_g = \partial V^4_g$ acts faithfully on the fundamental
group then the order of the group is bounded by a quadratic polynomial in $g$
(but not by a linear one), and hence the same holds also for 4-dimensional
handlebodies $V^4_g$.  As a consequence, each finite group $G$ acting on
$H^3_g$ or $V^4_g$ has a finite cyclic normal subgroup $G_0$ (the subgroup
acting trivially on the fundamental group) such that the order of $G/G_0$ is
bounded by a quadratic polynomial in $g$ ([Z4]).

\medskip

There arises naturally the question (as asked in [Z4]) whether there are
analogous polynomial bounds also for the orders of finite groups acting on
handlebodies $V^d_g$  of arbitrary dimension  $d$.  Whereas finite group
actions in dimension 3 are standard by the recent geometrization of such
actions after Thurston and Perelman, the situation in higher dimensions is more
complicated and not well-understood. Hence one is led to consider some kind of
standard actions also in higher dimensions. We will do so by uniformizing
handlebodies $V^d_g$ by Schottky groups (groups of M\"obius
transformations of the ball $B^d$ acting by isometries on its interior, the Poincar\'e-model
of hyperbolic space $\H^d$), thus realizing their interiors as complete
hyperbolic manifolds, and then considering finite groups of isometries of such hyperbolic
(Schottky) handlebodies  (see section 2 for the definition of Schottky groups).

\medskip

Our main results are as follows.

\bigskip

{\bf Theorem 1.}  {\sl  Let $G$ be a finite group of isometries of a
hyperbolic handlebody $V^d_g$ of dimension $d \ge 3$ and of genus $g > 1$ which acts
faithfully on the fundamental group. Then the order of $G$ is bounded by a polynomial of
degree $d/2$ in $g$ if $d$ is even, and of degree
$(d+1)/2$ if $d$ is odd. The degree $d/2$ is best possible in even dimensions whereas
in odd dimensions the optimal degree is at least $(d-1)/2$.}

\bigskip

By hypothesis such a group $G$ injects into the outer automorphism group of the fundamental
group of $V^d_g$, a free group of rank $g$. We note that by [WZ] the optimal upper bound
for the order of an arbitrary finite subgroup of the outer automorphism group Out$(F_g)$ of
a free group
$F_g$ of rank $g > 2$ is $2^g \, g!$ (i.e., exponential in $g$). It is shown in [Z2] that
every finite subgroup of Out$(F_g)$ can be induced (or realized in the sense of the Nielsen
realization problem) by an isomorphic group of isometries of a handlebody $V^d_g$ of
sufficiently high dimension $d$.

\medskip

Without the hypothesis that $G$ acts faithfully on the fundamental group, the
proof of Theorem 1 gives the following polynomial Jordan-type bound for finite
groups of isometries of $V^d_g$.

\bigskip

{\bf Corollary.}  {\sl Let $G$ be a finite group of isometries of a
hyperbolic handlebody $V^d_g$ of genus $g > 1$, and let $G_0$ denote the
normal subgroup of $G$ acting trivially on the fundamental group. Then the
following holds.

\medskip

i) $G_0$ is isomorphic to subgroup of the orthogonal group ${\rm SO}(d-2)$, and
the order of the factor group $G/G_0$ is bounded by a polynomial as in Theorem
1.

\medskip

ii) $G$ has a normal abelian subgroup, a subgroup of $G_0$, whose
index in $G$ is bounded by a polynomial as in Theorem 1.}

\bigskip

By the classical Jordan bound, each finite subgroup $G$ of a complex linear
group GL$(d,\Bbb C)$ has a normal abelian subgroup whose index in $G$ is
bounded  by a constant depending only on the dimension $d$ (see [C] for the optimal bound
for each $d$; see also [Z5] and its references for generalizations of
the Jordan bound in the context of diffeomorphism groups of manifolds).

\medskip

In more algebraic terms, Theorem 1 is equivalent to the following:

\bigskip

{\bf Theorem 2.}  {\sl Let $E$ be a group of M\"obius transformations of $S^{d-1}$ which is
a finite effective extension of a Schottky group ${\cal S}_g$ of rank $g > 1$. Then the
order of the factor group $E/{\cal S}_g$ is bounded by a polynomial in
$g$ as in Theorem 1.}

\bigskip

Here effective extension means that no element of $E$ acts trivially on ${\cal
S}_g$ by conjugation. By [Z2] every finite effective extension of a Schottky
group can be realized by a group of M\"obius transformations in some
sufficiently high dimension $d$.

\medskip

As a consequence of the geometrization of finite group actions in dimension
three, using the methods of [RZ, section
2] every finite group $G$ of diffeomorphisms of a 3-dimensional handlebody
$V^3_g$ is conjugate to a group of isometries, uniformizing $V^3_g$ by a
suitable Schottky group (which depends on $G$).  This is no longer true in higher
dimensions; however, if
$G$ is a finite group of diffeomorphisms of a 4-dimensional handlebody $V^4_g$ then,
uniformizing $V^4_g$ by a suitable Schottky group, $G$ acts also as a group of isometries of
$V^4_g$ inducing the same action on the fundamental group (applying the methods of
[Z4] to the boundary  3-manifold
$H^3_g$ of $V^4_g$). This raises naturally the following:

\bigskip

{\bf Questions.} i) Is every finite group $G$ of diffeomorphisms of a handlebody
$V^d_g$ isomorphic to a group of isometries of a hyperbolic handlebody $V^d_g$  (inducing
the same action on the fundamental group)?
 
\smallskip

ii) Is every finite group $G$ of diffeomorphisms of a ball
$B^d$ (i.e., a handlebody of genus zero) or of a sphere $S^{d-1}$ isomorphic to a subgroup
of the orthogonal group SO($d$)?

\medskip

In general, such a finite group $G$ of diffeomorphisms is not conjugate to a group of
isometries of a handlebody resp. to a group of orthgonal maps; we note that ii) is not true for
finite groups $G$ of homeomorphisms of $B^d$ or $S^{d-1}$, see [GMZ, section 7].

\bigskip

In section 2 we prove the first part of Theorem 1.  In section 3 we present examples of
finite isometric group actions on handlebodies which show that the degree $d/2$ of the
polynomial  bound in Theorem 1 is best possible in even dimensions (even for finite cyclic
groups $G$), and that a lower bound for the degree in  odd dimensions is
$(d-1)/2$. Note that for $d = 3$ the bound $(d+1)/2$ is not best possible  since it gives a
quadratic bound instead of the actual linear bound $12(g-1)$; for odd dimensions $d > 3$
we have no intuition at present if the optimal bound should be $(d-1)/2$ or $(d+1)/2$.

\bigskip

{\bf 2. Schottky groups and the Proof of Theorem 1.}
A {\it Schottky group} ${\cal S}_g$ of rank or genus $g$ is a group of
M\"obius transformations acting on a sphere $S^{d-1} = \partial B^d$ 
defined  in the following way  (analogously to the Schottky groups in dimension two acting
on $S^2$, see [L],[M] or [R, p. 584]; see also [Z2] for the following). Let $S_1, T_1,
\ldots, S_g, T_g$ be spheres of dimension $d-2$ in $S^{d-1}$ which bound disjoint balls
$B_1, D_1, \ldots, B_g, D_g$ of dimension $d-1$; choose M\"obius transformations $f_1,
\ldots, f_g$ such that
$f_i(S_i) = T_i$ and $f_i$ maps the exterior of $B_i$ to the interior of $D_i$. Then it is
easy to see that
$f_1, \ldots, f_g$ are free generators of a free group $\S_g$ of M\"obius
transformations.  The complement in $S^{d-1}$ of the interiors of the balls
$B_i$ and $D_i$ is a fundamental domain for the action of $\S_g$ on $S^{d-1} -
\Lambda({\cal S}_g)$ where $\Lambda({\cal S}_g)$ denotes the set of limit points of $\S_g$
in $S^{d-1}$ (a Cantor set). In this definition, one may consider round spheres $S_1, T_1,
\ldots, S_g, T_g$ (thus defining a so-called classical Schottky group), or just
topological spheres (and it is known that non-classical Schottky groups esist); however this
is not relevant for the present paper, in particular in the examples constructed in section
3 the Schottky subgroups will be always classical).

\medskip

The group of M\"obius transformations of $S^{d-1}$ extends naturally to the
interior of the ball $B^d$ ("Poincar\'e extension") where it becomes the group
of orientation-preserving isometries of the Poincar\'e-model of hyperbolic
space $\Bbb H^d$. The action of ${\cal S}_g$ is free and properly discontinuous
on the interior $\Bbb H^d$ of $B^d$, and a fundamental domain for this action
is the region of $\Bbb H^d$ bounded by all hyperbolic hyperplanes defined by
the spheres $S_i$ and $T_i$ (i.e., half-spheres of dimension $d-1$ orthogonal
to $S^{d-1}$ along these spheres).  The quotient $(B^d - \Lambda({\cal
S}_g))/{\cal S}_g$ is a handlebody $V^d_g$ whose interior $\Bbb H^d/{\cal S}_g$
has the structure of a complete hyperbolic manifold, and we say that the
Schottky group ${\cal S}_g$ uniformizes the handlebody $V^d_g$. When speaking
of a finite group $G$ of isometries of a handlebody $V^d_g$ we then intend that
$V^d_g$ can be uniformized by a Schottky group ${\cal S}_g$ such that $G$ acts
by hyperbolic isometries on the interior of $V^d_g$.

\medskip

Let $V^d_g$ be a handlebody uniformized by a Schottky group $\S_g$. Let $G$ be
a finite group of isometries of $V^d_g$ which induces a faithful action on the
fundamental group. The group of all lifts of elements of $G$ to the universal
covering $B^d - \Lambda({\cal S}_g)$ of $V^d_g$ defines a group $E$ of M\"obius
transformations of $B^d$, with factor group $E/\S_g \cong  G$, so we have a
finite extension
$$1 \to \S_g  \hookrightarrow  E  \to G \to 1;$$
by general covering space theory, this extension is effective since $G$ acts
faithfully on the fundamental group of $V^d_g$ (isomorphic to the group $\S_g$
of covering transformations).

\bigskip

{\bf Lemma 1.}  {\sl The extension $1 \to \S_g  \hookrightarrow  E  \to G \to 1$ is
effective if and only if $E$ has no non-trivial finite normal subgroups.}

\medskip

{\it Proof.}  Let $F$ be a finite normal subgroup $E$. Since the
intersection of $F$ with the normal torsionfree subgroup $\S_g$ of $E$ is trivial, the
normal subgroups $F$ and $\S_g$ of $E$ commute elementwise (any commutator $fsf^{-1}s^{-1}$
of elements $f \in F$ and $s \in \S_g$ is an element of both $F$ and $\S_g$ and hence
trivial). Hence if if the extension is effective, $F$ has to be trivial.

\medskip

Conversely, suppose that every finite normal subgroup of $E$ is trivial. The
subgroup of elements of the finite extension $E$ of $\S_g$ inducing by
conjugation the trivial automorphism of $\S_g$ is clearly finite (since the
center of $\S_g$ is trivial), normal and hence trivial, so the extension is
effective.

\medskip

This completes the proof of Lemma 1.

\bigskip

As a consequence of Stalling's structure theorem for groups with infinitely many ends, a
finite extension $E$ of a free group is  the fundamental group $\pi_1\G$ of a
finite graph of finite groups  $\G$ ([KPS]); here $\Gamma$ denotes a finite graph, and to
its vertices $v$ and edges $e$ are associated finite vertex groups $G_v$ and edge groups
$G_e$, with inclusions of the edge groups into the adjacent vertex groups. The fundamental
group $\pi_1\G$ of the finite graph of finite groups $\G$ is the iterated free product with
amalgamation and HNN-extension of the vertex groups amalgamated over the edge groups, first
taking the iterated free product with amalgamation over a maximal tree of $\Gamma$ and then
associating an HNN-generator to each of the remaining edges. We note that each finite
subgroup of $E =
\pi_1\G$ is conjugate into a vertex group of $\G$, and that the vertex groups are maximal
finite subgroups of $E$  (see [ScW], [Se] or [Z3] for the standard theory of graphs of
groups and their fundamental groups).

\medskip

We will assume in the following that the graph of
groups $\G$ has no {\it trivial edges}, i.e. no edges with two different vertices
such that the edge group coincides with one of the two vertex groups (by
collapsing trivial edges, i.e. amalgamating the two vertex groups into a single vertex
group); we say that such a graph of groups is in {\it normal form}.

\medskip

We denote by
$$\chi\G = \sum {1 \over |G_v|} - \sum {1 \over |G_e|}$$
the  {\it Euler characteristic} of the graph of groups $\G$ (the sum is taken over all
vertex groups $G_v$ resp. edge groups $G_e$ of $\G$); then, by multiplicativity of Euler
characteristics under finite coverings of graphs of groups,
$$g-1 =  -\chi\G \; |G|$$
(see [ScW] or [Z3]); note that this is positive since we are assuming that $g > 1$.

\medskip

The finite extension $E = \pi_1\G$ of the Schottky group $\S_g$ is a group of
M\"obius transformations of $B^d$ and acts as a group of hyperbolic isometries
on its interior $\H^d$. Each finite group of isometries of hyperbolic space
$\H^d$ has a global fixed point in $\H^d$ and is conjugate to a finite group of
orthogonal transformations of $B^d$ (which are exactly the isometries of $\H^d$
which fix the origin in $B^d$).  In particular each finite vertex group $G_v$
of $E = \pi_1\G$ has a fixed point in $\H^d$ and is isomorphic (conjugate) to a
subgroup of the orthogonal group O($d$), and different vertex groups of $\G$
have different fixed points (since the vertex groups are maximal finite
subgroups of $E$ and the action of $E$ is properly discontinuous in $\Bbb
H^d$); also, if a vertex group fixes a point in $\H^d$ then it is the maximal
finite subgroup of $E$ fixing this point.

\medskip

Consider a non-closed edge $e$ of $\G$, i.e. with two distinct vertices $v_1$ and $v_2$,
with edge group $G_e$ and vertex groups $G_1$ and $G_2$ (which we consider as subgroups of
$E$), with  $G_e = G_1 \cap G_2$.  Let $P_1 \ne P_2$ be fixed points of $G_1$ resp.
$G_2$ in $\H^d$; then $P_1$ and $P_2$ define a
hyperbolic line $L$ which is fixed pointwise by the edge group $G_e = G_1 \cap G_2$. The
line $L$ intersects $S^{d-1} = \partial B^d$ in two points which are fixed by $G_e$;
moreover no subgroup of $G_1$ larger than $G_e$ can fix one of these two points
since otherwise it would fix pointwise the line $L$ and hence $P_2$, so it would be
contained also in $G_2$.

\medskip

Now let $e$ be a closed edge of $\G$, i.e. an edge with only one vertex $v$. There are two
inclusions of the edge group $G_e$ into the vertex group $G_v$ defining two subgroups
$G_e$ and $G_e'$ of $G_v$; denoting by $t$ an HNN-generator
corresponding to the edge $e$, we have that $t^{-1}G_e't = G_e$, and
$G_e = G_v \cap (t^{-1}G_vt)$.  Note that $t$ has
infinite order so it does not fix any point in $\H^d$. Let $P$ be a fixed point of the
finite subgroup $G_v$ of $E$ in $\H^d$; then $t^{-1}G_vt$ fixes the
point $t(P) \ne P$, and  its subgroup $G_e = t^{-1}G_e't$ fixes the hyperbolic line $L$
defined by $P$ and $t(P)$. As before, the hyperbolic line $L$ intersects
$S^{d-1} = \partial B^d$ in two points which are fixed by $G_e$, and $G_e$ is the maximal
subgroup of $G_v$ fixing these two points.

\medskip

Note also that, since $G_e$ fixes a point in $S^{d-1}$, it is in fact
isomorphic (conjugate) to a subgroup of the orthogonal group SO($d-1$).
Summarizing, we have:

\bigskip

{\bf Lemma 2.}  {\sl Let $G_v \subset E$ be a vertex group of the graph of groups $\G$, and
let $G_e \subset G_v$ be an adjacent edge group. Then $G_v$ has a global fixed point in
$\Bbb H^d$, and $G_e$ has a global fixed point in
$S^{d-1} = \partial B^d$ which is not fixed by any other element of $G_v$. In particular,
every vertex group is isomorphic to a subgroup of the orthogonal group ${\rm SO}(d)$, and
every edge group is isomorphic to a subgroup of ${\rm SO}(d-1)$.}

\bigskip

We need also the following lemma which is contained in [Z4, proof of Theorem
1]; since its proof is short, we present it for the convenience of the reader.
Let $\chi = \chi\G$ denote the Euler characteristic of $\G$; note that $-\chi >
0$ since $g>1$, and that for any graph of groups in normal form one has $-\chi
\ge 0$ unless it consists of a single vertex.

\bigskip

{\bf Lemma 3.}  {\sl  Let $e$ be an edge of $\Gamma$. Denote by $n$ the order of $G$ and
by $a$ the order of the edge group $G_e$. Then
$${n \over a} \le 6(g-1).$$}

\medskip

{\it Proof.}  Suppose first that $e$ is a closed edge. If $e$ is the only edge of
$\G$ then
$$-\chi  \ge  {1 \over a} - {1 \over 2a} = {1 \over a}, \hskip 5mm
g-1 = -\chi n \ge {n \over 2a},  \hskip 5mm  {n \over a} \le 2(g-1).$$
If $e$ is closed and not the only edge then
$$-\chi  \ge  {1 \over a}, \hskip 5mm
g-1 = -\chi n \ge {n \over a},  \hskip 5mm  {n \over a} \le g-1.$$

Suppose that $e$ is not closed. If $e$ is the only edge of $\G$ then both
vertices of $e$ are isolated and
$$-\chi  \ge {1 \over a} - {1 \over 2a} - {1 \over 3a} = {1 \over 6a}, \hskip 5mm   g-1 =
-\chi \; n  \; \ge  \; {n \over 6a}, \hskip 5mm  {n \over a} \le 6(g-1).$$

\medskip

If $e$ is not closed, not the only edge and has exactly one isolated vertex then
$$-\chi  \ge {1 \over a} - {1 \over 2a} = {1 \over 2a}, \hskip 5mm   g-1 =
-\chi \; n  \; \ge  \; {n \over 2a}, \hskip 5mm  {n \over a} \le 2(g-1).$$
Finally, if $e$ is not closed, not the only edge and has no isolated vertex then
$$-\chi  \ge {1 \over a}, \hskip 5mm   g-1 =
-\chi \; n  \; \ge  \; {n \over a}, \hskip 5mm  {n \over a} \le g-1.$$
Concluding, in all cases the inequality of Lemma 3 holds, completing the proof of the
lemma.

\bigskip

After these preparations, we can now start with the  actual:

\bigskip

{\it Proof of Theorem 1.}  Let $e$ be any edge of the finite graph of finite
groups $\G$ given by the $G$-action.  By Lemma 2, $G_e$ has a global fixed
point in $S^{d-1} = \partial B^d$ and is isomorphic to a subgroup of the
orthogonal group SO($d-1$). By the classical Jordan bound for subgroups of
GL$(d-1,\Bbb C)$, the edge group $G_e$ has an abelian subgroup $A_1$ whose
index in $G_e$ is bounded by a constant $c$ depending only on the dimension. We
will find a polynomial upper bound in $g$ for the order $a_1$ of the abelian
group $A_1$; this will imply then a polynomial bound of the same degree also
for the order $a \le c \, a_1$ of $G_e$, and finally for the order $n$ of $G$
since, by Lemma 3,
$$n \;  \le \; 6(g-1) \, a \; \le \; c\, 6(g-1)\, a_1.$$

\medskip

Let $E_1$ be the subgroup of $E$ generated by $\S_g$ and $A_1$ (which is again an effective
extension of $\S_g$, with factor group $A_1$). Then also $E_1$ is the fundamental group of a
finite graph of finite groups in normal form which we denote again by $\G$. Since the finite
group $A_1$ has a fixed point in $\H^d$, up to conjugation it is the vertex
groups $G_v$ of some vertex $v$ of $\G$, and its fixed point set in $S^{d-1}$ is a sphere
$S^{d_1}$ of dimension $d_1 \ge 0$ (since  $G_e$ has a global fixed point in
$S^{d-1}$). Since  $\G$ has no trivial edges and $E_1$ has no non-trivial finite normal
subgroups by Lemma 1, some edge adjacent to
$v$ has an edge group $A_2$ of order $a_2 < a_1$ (i.e., properly contained in
$A_1$). By Lemma 3,
$$a_1 \; \le \; 6(g-1) \, a_2.$$
By Lemma 2, the edge group $A_2$ has a fixed point in $S^{d-1} = \partial B^d$ which is not
fixed by any other element of the vertex group $A_1$, hence the fixed point set of $A_2$ in
$S^{d-1}$ is a sphere $S^{d_2}$ of dimension  $d_2 > d_1$.

\medskip

We iterate the construction and consider the subgroup $E_2$ of $E_1$ generated by $\S_g$
and $A_2$, obtaining an edge group $A_3$ for $E_2$ which fixes a sphere $S^{d_3}$ of
dimension $d_3 > d_2$ in $S^{d-1}$, of order
$$a_2 \; \le \; 6(g-1) \, a_3.$$
Hence, after at most $d-1$ steps, we end up with a trivial edge
group fixing all of $S^{d-1}$. Collecting, we obtain the polynomial bound
$$n \; \le \; c \, 6^d(g-1)^d$$
of degree $d$ in $g$ for the order of $G$.

\bigskip

In order to obtain a polynomial bound of the degree given in Theorem 1 we argue as
follows. Suppose  that the fixed point set of the normal subgroup
$A_2$ of $A_1$ is a sphere $S^{d_1+1}$ of dimension $d_2 = d_1+1$; note that
$S^{d_1+1}$ is invariant under the action of $A_1$. Let
$A_1'$ denote the subgroup of index one or two of $A_1$ which acts orientation-preservingly
on $S^{d_1+1}$.  Then  $A_1'$ fixes  $S^{d_1+1}$ pointwise since
otherwise the fixed point set of $A_1'$ would be a sphere of
codimension at least two in $S^{d_1+1}$; this is not possible since already
$A_1$ has fixed point set $S^{d_1}$ of dimension $d_1$.  Continuing now with $A_1'$ in the
place of $A_1$, we can assume that the dimensions $d_i$ increase by at least two in each
step.  Hence the number of steps is at most $d/2$ if $d$ is even, and $(d+1)/2$ if $d$
is odd, and  this gives the degree of the polynomial upper bound as stated in Theorem 1.

\medskip

This completes the proof of the first part of Theorem 1; the second part on the optimality
of the degree $d/2$ for even $g$ and the lower bound $(d-1)/2$ for odd $g$ will follow from
the examples of finite group actions on handlebodies constructed in the next section.

\bigskip

{\it Proof of the Corollary.}  The proof proceeds along the lines of the proof
of Theorem 1, with the following difference.  In the proof of Theorem 1 we
consider the sequence of abelian subgroups $A_1, A_2, \ldotsÊ\;$ of $G$; after
finitely many steps, this ended with the trivial group, using the effectiveness
of the corresponding extensions $E_1, E_2, \ldots \;$ of $\S_g$. Without
effectiveness, the sequence $A_1, A_2, \ldots \;$ of $G$ ends with an abelian
group $A_m$ which is a normal subgroup of the corresponding extension $E_m$; in
particular, $A_m$ acts trivially on $\S_g$ and is a subgroup of $G_0$. The
index of $A_m$ in $G$ is bounded by a polynomial as in the proof of Theorem 1,
hence also the index of $G_0$ in $G$ is bounded by such a polynomial.

\medskip

The group $G_0$ lifts to an isomorphic normal subgroup of the
extension $E$ of $\S_g$ which we denote also by $G_0$. The finite group $G_0$ has a fixed
point in $\Bbb H^d$; we can assume that it fixed the origin $O \in B^d$ and hence is
isomorphic to a subgroup of SO($d$). Since $G_0$ is normal in $E$, it is contained (up to
conjugation) in each edge group of the graph of groups $\G$.
By Lemma 2, $G_0$ has a global fixed point also in $S^{d-1} = \partial B^d$, hence it fixes
pointwise a great sphere of dimension at least zero in $S^{d-1}$, and a linear
subspace $B$ of dimension at least one in $B^d$. Since $G_0$ commutes elementwise
with $\S_g$, the Schottky group $\S_g$ acts on $B$. Since the action of $\S_g$ is
properly discontinuous and $g>1$, $B$ has dimension at least two. Now $G_0$ acts also on
the orthogonal complement of $B$ in $O \in B^d$, a linear subspace of codimension at
least two, so $G_0$ is isomorphic to a subgroup of the orthogonal group SO($d-2$).

\medskip

Finally, by the classical Jordan bound for linear groups, the subgroup $G_0$ of SO($d-2$)
contains a normal abelian subgroup whose index is bounded by a constant depending only on
the dimension $d$. By taking the intersection of this normal abelian subgroup with all
isomorphic normal subgroups of $G_0$ we obtain a characteristic abelian
subgroup $A$ of $G_0$ whose index in $G_0$ is also bounded by a constant depending only on
the dimension $d$. Hence the indices of $A$ and $G_0$ in $G$ are bounded by polynomials in
$g$ of the same degree.

\medskip

This completes the proof of the Corollary.

\bigskip

{\bf 3. Examples.}  We construct isometric actions of finite groups $G$ on
handlebodies which realize the lower bounds for the degrees of the polynomial bounds in
Theorem 1; specifically, we prove the following:

\bigskip

{\bf Proposition.}  {\sl  For a fixed $k \ge 2$ and all $m \ge 2$, the finite group
$G = (\Bbb Z_m)^k$  admits an action, faithful on the fundamental group, on a
handlebody
$V^d_g$ of genus $g = mk-k$ and dimension $d = 2k$ and $2k+1$; in particular, the
order $n = m^k$ of $G$ is given by the polynomial
$$ n \; = \; (g + k)^k/k^k \; = \; (1 +  g/k)^k$$
of degree $k = d/2$ in $g$ if $d$ is even, and $k = (d-1)/2$ if $d$
is odd.}

\bigskip

{\it Proof.}  For $k > 1$, let $G = C_1 \times \ldots \times C_k  \cong (\Bbb Z_m)^k$, of
order $n = m^k$, be the product of $k$ cyclic groups $C_i \cong \Bbb Z_m$ of
order $m$. Choose an orthogonal action of $G$ on the closed ball $B^{2k}
\subset \Bbb R^{2k}$ of dimension $d = 2k$ as follows. Decomposing $\Bbb R^{2k}
= P_1 \times \ldots \times P_k$ as the product of $k$ orthogonal planes $P_i$,
each $C_i$ acts on $P_i$ faithfully by rotations and trivially on the $k-1$
orthogonal planes.

\medskip

Define a finite graph of finite groups $\G$ as follows.  The graph $\Gamma$ is a
star-shaped graph with one central vertex $v$ with vertex group
$G_v = G = C_1 \times \ldots \times C_k$ and $k$ non-closed edges  $e_1, \ldots, e_k$ each
having $v$ as a vertex, with edge groups
$$G_{e_1} = C_2 \times \ldots \times C_k \;, \hskip 2mm  G_{e_2} = C_1 \times C_3 \times
\ldots \times C_k \;, \hskip 2mm  \ldots \;,  \hskip 2mm G_{e_k} = C_1 \times \ldots \times
C_{k-1}$$ (i.e., exactly $C_i$ is missing in $G_{e_i}$). Hence $\Gamma$ has $k+1$
vertices, by definition all with vertex group $G = C_1 \times \ldots \times C_k$, and the
Euler characteristic of $\G$ is
$$ \chi \; = \; (k+1){1 \over m^k} \,  - \,  k{1 \over m^{k-1}} \;.$$

There is an obvious projection of the fundamental group $E = \pi_1\G$ of the graph of groups
$\G$ onto $G$; its kernel is a free group $F_g$ of some rank $g$, and we have an extension
$$ 1 \to F_g  \hookrightarrow  E  \to G  \to  1$$
which by construction of $\G$ is effective (has no nontrivial finite normal subgroups, see
Lemma 1). The rank $g$ is given by
$$g - 1 \, = \, (-\chi) n \, = \, (-\chi) m^k \, = \, mk - (k+1), \hskip 3mm g \, =
\, mk-k,$$ hence
$$ n \, = \, m^k \, = \, (g + k)^k/k^k$$
which is a polynomial of degree $k = d/2$ in $g$ and gives the maximal
possibility for the degree in Theorem 1 for even dimensions $d$.

\medskip

We realize $E = \pi_1\G$ as a group of M\"obius transformations of $B^d$, $d =
2k$, such that its subgroup $F_g$ corresponds to a Schottky group $\S_g$. Then the quotient
$(B^d -  \Lambda({\cal S}_g))/{\cal S}_g$ is a handlebody $V_g^d$ of genus $g$, and $E$
projects to an action of the factor group
$E/\S_g \cong G$ on $V_g^d$ which is faithful on the fundamental group. In particular,
the degree $d/2$ in Theorem 1 is best possible for even dimensions $d = 2k$.

\medskip

The realization of $E = \pi_1\G$ as a group of M\"obius transformations of $B^d$
proceeds inductively  by standard combination methods (similar as in [Z2, section 3]).
Starting with the orthogonal group $G$ described above, we realize first the
free product with amalgamation
$$G_v *_{G_{e_1}}G_{v_1} = G *_{G_e}G_1$$
where $e = e_1$ denotes the first edge of $\Gamma$, with vertices $v$ and $v_1$
and vertex groups $G = G_v$ and $G_1 = G_{v_1} \cong G$. By construction, the
fixed point set of the subgroup $G_e$ of $G$ is a 2-ball $B_1$ in $B^d$
defining a hyperbolic plane in $\Bbb H^d$ which will be denoted also by $B_1$. Let $L_1$ be
a hyperbolic half-line in $B_1$ starting from its center 0 and ending in a point  $R_1$ in
$S^{d-1} =
\partial B^d$. Let $V_1$ be a neighbourhood of
$R_1$ in $B^d$ bounded by a hyperbolic hyperplane $H_1$ in $\Bbb H^d$
orthogonal to $L_1$; choose $V_1$ sufficiently small such that  $f(V)$ is
disjoint form $V$ for all $f \in G - G_e$ (note that $G_e$ fixes $L_1$
pointwise but that no larger subgroup of $G$ fixes $L_1$ by construction of
$G$). The reflection $\tau_1$ in the  hyperbolic hyperplane $H_1$ commutes
elementwise with $G_e \subset G$ and, considering  $G_1 = \tau_1 G
\tau_1^{-1}$, we have that  $G \cap G_1 = G_e$. Similar as for Schottky groups
it is now easy to see that the group of M\"obius transformations generated by
$G$ and $G_1$ is isomorphic to the free product with amalgamation $G
*_{G_e}G_1$, and that every torsionfree subgroup of finite index is in fact a
Schottky group (cf. [Z2] and the combination theorems in [M]).

\medskip

We iterate the construction and adjoin $G_{e_2}$. Let $L_2$ be a hyperbolic half-line
starting in the center 0 and ending in a point $R_2$ of  $S^{d-1} = \partial B^d$ such that
$R_2$ does not lie in $G(V_1)$. Let $V_2$ be a small neighbourhood of $R_2$ in $B^d$,
bounded by a hyperbolic hyperplane $H_2$ orthogonal to $L_2$
which does not intersect $G(V_1)$. With $G_2 =  \tau_2 G \tau_2^{-1}$ where
$\tau_2$ denotes the reflection in $H_2$, this realizes the free product with amalgamation
$$G_{v_2} *_{G_{e_2}}G_v *_{G_{e_1}}G_{v_1}$$
 as a group of M\"obius transformations.
Continuing in this way, after
$k$ steps $E$ is realized as a group of M\"obius transformations, with $F_g$
corresponding to a Schottky group $\S_g$.

\bigskip

Finally, in odd dimensions $d = 2k+1$, we extend the orthogonal action of $G$
on $B^{2k}$ described above to an orthogonal action on $B^{2k+1}$ (trivial on the last
coordinate) and then proceed as before. We get a polynomial of degree $k = (d-1)/2$ in $g$
for the order $n$ of $G$ whereas Theorem 1 gives a polynomial bound of degree
$(d+1)/2$. As noted in the introduction, the optimal degree in dimension $d =
3$ is in fact 1, but for odd dimensions  $d > 3$ it remains open.

\medskip

This completes the proof of the Proposition, and also of Theorem 1.

\bigskip

The examples given in the Proposition are for finite abelian groups $G$. By
suitably modifying the construction, one obtains also examples for finite
cyclic groups as follows.

\medskip

Let $d = 2k$ be a fixed even dimension, and let $p > k$ be any prime. For $i =
1, \ldots,k$, the  $k$ integers $q_i = p + i\; k!$, , are pairwise coprime: in
fact, if a prime $p'$ divides $q_i$ then $p' > k$; if $p'$ divides also $q_j$,
for some  $j > i$, then $p'$ divides $q_j - q_i =  (j-i)\;  k!$ which is a
contradiction. Then $G = \Z_{q_1} \times \ldots \times \Z_{q_k}$ is a cyclic
group of order $n = q_1 \ldots q_k$. In analogy with the proof of the
Proposition, let $\G$ be a star-shaped graph of groups with $k+1$ vertices all
with vertex group $G$, and with $k$ edges where in each edge group is missing
exactly one of the factors $\Z_{q_i}$ of $G$, with

$$\chi = \chi\G \; = \; {k+1 \over n} -  {q_1 \over n}- \ldots - {q_k \over n}.$$
There is an obvious surjection of $\pi_1\G$ onto $G$, its kernel is a free group of rank
$g$ with
$$g - 1 \; = \;  (-\chi) \, n  \; =  \;  -(k+1) + q_1 + \ldots + q_k,$$
$$g \; = \; -k + kp + (1+ \ldots + k) \, k!,$$
$$ p = (g + c_k)/ k),$$
for a constant $c_k$ depending only on $k$. Now
$$|G| \; = \; n \; = \; q_1 \ldots q_k \; \ge \; p^k \; \ge \; (g+c_k)^ k/k^k,$$
so the order of $G$ is bounded from below by a polynomial of degree $k = d/2$ in $g$.

\medskip

Finally, the geometric realizations of $G$ and $E = \pi_1\G$ are exactly as in the proof of
the Proposition.

\bigskip \bigskip

{\bf Acknowledgment.}  The authors were supported by a FRV grant from
Universit\`{a} degli Studi di Trieste.

\bigskip \bigskip

\centerline {\bf References}

\bigskip

\item {[C]}  M.J. Collins, {\it  On Jordan's theorem for complex linear groups,}
J. Group Theory 10   (2007),  411-423

\item {[GMZ]} A. Guazzi, M. Mecchia, B. Zimmermann, {\it On finite groups acting on
acyclic low-dimensional manifolds,}  Fund. Math. 215  (2011),  203-217

\item {[KPS]}  A. Karras, A. Pietrowski, D. Solitar, {\it  Finite and infinite extensions
of free groups,} J. Austral. Math. Soc. 16  (1972),  458-466

\item {[L]}  J. Lehner, {\it  Discontinuous Groups and Automorphic Functions,} Math.
Surveys 8, Amer. Math. Soc. 1964

\item {[M]}  B. Maskit, {\it  Kleinian groups,} Grundlehren Math. Wiss. 287,
Springer-Verlag  1988

\item {[MMZ]} D. McCullough, A. Miller, B. Zimmermann,  {\it Group actions on
handlebodies,}  Proc. London Math. Soc.  59   (1989), 373-415

\item {[MZ]} A. Miller, B. Zimmermann,  {\it  Large groups of symmetries of
handlebodies,}  Proc. Amer. Math. Soc. 106  (1989),  829-838

\item {[R]} J.G. Ratcliffe,  {\it  Foundations of Hyperbolic Manifolds,}
Graduate Texts in Mathematics 149, Springer-Verlag 1980

\item {[RZ]} M. Reni, B. Zimmermann,  {\it Handlebody orbifolds and Schottky
uniformization of hyperbolic 2-orbifolds,}  Proc. Amer. Math. Soc. 123 (1995),
3907-3914

\item {[ScW]} P. Scott, T. Wall,  {\it  Topological methods in group theory,}
Homological Group Theory, London Math. Soc. Lecture Notes 36, Cambridge University
Press  1979

\item {[Se]} J.P. Serre,  {\it  Trees,}  Springer 1980

\item {[WZ]} S. Wang, B. Zimmermann,  {\it  The maximum order finite groups of
outer automorphisms of free groups,}  Math. Z. 216  (1994), 83-87

\item {[Z1]} B. Zimmermann,  {\it \"Uber Abbildungsklassen von
Henkelk\"orpern,}  Arch. Math. 33  (1979),  379-382

\item {[Z2]} B. Zimmermann,  {\it \"Uber Hom\"oomorphismen n-dimensionaler Henkelk\"orper
und endliche Erweiterungen von Schottky-Gruppen,}  Comm. Math. Helv. 56   (1981), 474-486

\item {[Z3]} B. Zimmermann,  {\it  Generators and relations for discontinuous groups,}
Generators and Relations in Groups and Geometries (eds. Barlotti, Ellers, Plaumann,
Strambach),  NATO Advanced Study Institute Series 333,  Kluwer Academic
Publishers  1991,   407-436

\item {[Z4]}   B. Zimmermann,  {\it  On finite groups acting on a connected sum of
3-manifolds $S^2 \times S^1$,} Fund. Math. 226 (2014), 131-142

\item {[Z5]}   B. Zimmermann,  {\it  On Jordan type bounds for finite groups acting on
compact 3-mani-folds,} Arch. Math. 103 (2014), 195-200

\bigskip

Dipartimento di Matematica e Geoscienze

Universit\`a degli Studi di Trieste

34127 Trieste, Italy

E-mail:  mmecchia@units.it,  zimmer@units.it

\bye